\documentclass[12pt]{article}
\usepackage{amsmath}
\usepackage{amssymb}
\usepackage{amsfonts}
\usepackage{eucal}
\numberwithin{equation}{section}
\newtheorem{Theorem}{Theorem}[section]
\newtheorem{Proposition}[Theorem]{Proposition}
\newtheorem{Definition}[Theorem]{Definition}

\newtheorem{Remark}[Theorem]{Remark}
\begin{document}
\title{A randomized first-passage problem for drifted Brownian motion subject to hold and jump  from a boundary}
\author{Mario Abundo\thanks{Dipartimento di Matematica, Universit\`a  ``Tor Vergata'', via della Ricerca Scientifica, I-00133 Rome, Italy.
E-mail: \tt{abundo@mat.uniroma2.it}}
}
\date{}
\maketitle

\begin{abstract}
\noindent We study an inverse first-passage-time problem for
Wiener process $X(t)$ subject to hold and jump from a boundary
$c.$ Let be given a threshold  $S>X(0) \ge c,$ and a distribution
function $F$ on $[0, + \infty ).$ The problem consists in finding
the distribution of the holding time at $c$ and the distribution
of jumps from $c,$ so that the first-passage time of $X(t)$
through $S$ has distribution $F.$
\end{abstract}


\section{Introduction}
In this brief note, we generalize the inverse
first-passage-time problem for Wiener process $X$ with random jumps from a boundary $c$,
studied in our paper \cite{abundo:saa13}, by introducing holding when $X$ hits $c.$
Diffusions with holding and jumping boundary (DHJ) were first introduced by Feller (\cite{feller:1954}), and their
ergodic properties were studied in \cite{peng:2013}, and in \cite{benari:spa08} in the case without  holding.
In the one dimensional case, a DHJ coincides with a
diffusion in $I=(c,d),$ until it hits the boundary $c;$ when this occurs, the process waits there for an
independent random duration following an exponential distribution with parameter $\beta >0,$ and then it makes a random jump in the
interior of $I,$ according to a given distribution $\nu, $ and starts the diffusion afresh.
The (direct) first-passage time (FPT) of a one dimensional DHJ through a barrier $S$ with  $c \le X(0) < S \le d ,$ was investigated in
\cite{peng:2014},
while in \cite{abundo:saa11} the FPT was studied in the case without holding, i.e. $\beta = + \infty $ (see also \cite{peng:sta09}).
\par
The drifted Brownian motion with holding and jump boundary $c$ can
be described as follows. Let $\widetilde X _ \mu (t) = x+ \mu t +
B_t ,$ where $B$ is standard  Brownian motion (BM), and $x \in I =
[c,d].$ We construct a new stochastic process $ X,$ called drifted
Brownian motion with holding and jump boundary $c$: starting from
$x \in I, $ at time $t=0,$ we set $X  (t)= \widetilde X _\mu (t)$
until the (random) time at which the boundary $c$ is reached, then
the process $\widetilde X _ \mu $ is killed and it is continued as
the new process $X  $ which holds for an independent exponential
time at $c ,$ with the holding rate $\beta ,$ then it makes a
random jump from  $c$ inside $I$ and it starts afresh
(independently of the past history) from a point $U \in I$ which
is a random variable with distribution $\nu,$ and then it evolves
as  $\widetilde X _ \mu ,$ until the boundary $c$ is reached
again, and so on. Indeed, in the present paper we take $c=0$ and
$d>0,$  for the sake of simplicity. \par
Hitting the boundary
$c=0$ by the process $X$ can be interpreted as the occurrence of a
catastrophe with delayed effect; after holding, it consists in
resetting the state of $ X$ to a new value, and this has various
applications in the scope of diffusion models in Queueing theory
(see e.g. \cite{dicresc:que03}), as well as in Mathematical
Finance  and in Biology (see \cite{abundo:saa11}). As an example,
$X(t)$ can represent the number of individuals of a population at
time $t$ (for instance fishes in a little lake); whenever the
population goes extinct (i.e. it reaches the level $0$), after a
waiting time, it is performed a restocking by  a random amount
$U>0.$ \par\noindent Notice that, considering $\widetilde X _ \mu
(0) \ge c  $ and  $c=0 ,$ is not restrictive. For instance, let us
consider the diffusion $Z$ with holding and jumping boundary $c
\ge Z(0),$ associated to $\widetilde Z(t)= \widetilde Z(0) + \mu t
+ B_t$ with  $\mu >0, \ Z(0) >0 ,$ and $U \in (0,c). \ Z$ is similar to the
Wiener-type neuronal model in the presence of refractoriness,
considered in \cite{alb:2007}, \cite{peng:2014}, for which the
neuron fires when its voltage exceeds the threshold $c$, and after
the refractoriness period, the voltage is reset to $U \in (0,c).$
We observe that the process $Z$ can be reduced to our case;
indeed, by using that $-B_t$ is distributed as $B_t,$ it follows
that the distribution of the first hitting time of $\widetilde Z$
to $c,$ when starting from $\widetilde Z(0) \le c,$ is nothing but
the distribution of the first hitting time of $\widetilde X _ {-
\mu}$ to zero, when starting from $c-\widetilde Z(0) \ge 0.$
 \par
Let  $S \in I=(0, d)$ be a constant barrier, and suppose that the
holding rate $\beta$  and  the jump distribution $\nu (u) = P( U
\le u), \ u \in (0,S)$ are assigned;  then, the FPT of $X$ over
$S,$ when starting from $0< x <S,$ is:
\begin{equation} \label{tausuno}
\tau_S(x)=  \tau _ {S, \beta, \nu} (x) = \inf \{ t >0 :  X (t) > S | X(0)= \widetilde X(0)=x \} .
\end{equation}
The direct FPT problem for $ X$ was studied in \cite{peng:2014}, while the  FPT for
the process without holding (i.e. $\beta = + \infty)$ was studied in
 \cite{abundo:saa11}; here, we deal with the following inverse FPT (IFPT) problem:
\bigskip

{\it for a given barrier $S \in I $ and  $0 < x< S,$ let  $F$ be a  distribution on $[0, + \infty ),$ \par
then the IFPT problem  consists in finding the holding rate $\beta $ and the jump distribution $\nu$ \par
or its density $g,$ if it  exists, so  that
$\tau _ S (x)$ has distribution $F.$ } \bigskip

\noindent
In this case, we will say that the the pair $(g, \beta)$ is solution to the IFPT problem; notice that the density $g$ has support $(0,S).$ \par\noindent
This IFPT problem is analogous to that considered in \cite{abundo:saa13}
 for diffusion with only jumps (without holding) from a boundary,
and studied in  \cite{abundo:stapro13}, \cite{abundo:stapro12}, \cite{jackson:stapro09}
for diffusions without holding and without jumps, and in \cite{abundo:saa14} for reflected diffusions.
It has interesting applications in
Mathematical Finance, in particular in credit risk modeling, where
the FPT represents a default event of an obligor (see \cite{jackson:stapro09}), and in Biology  (see \cite{lansky:bio89}).
The paper is organized
as follows: Section 2  contains  the main results; Section 3 provides some extension to more general diffusions and in
 Section 4 we report some explicit examples.

\section{Main results}
Let $\widetilde X $ be a one-dimensional diffusion on $D= (0, + \infty)$ and let
$X$ the associated diffusion with holding and jumps from the boundary $c=0,$ starting from $x >0.$
Following the notations of \cite{peng:2014},
we denote by $\beta >0$ the holding rate of $X$ at $0,$ and by $\nu$ the jump distribution when it hits $0.$
Moreover, let $\tau (x)$ be the first-passage time (FPT) of $X$ through $S,$ with $x < S,$ and let $h_S (t,x)$ the density of $\tau (x)$
and  $ \widehat h _S ( \lambda, x )$ its Laplace transform $( \lambda >0).$ Finally, denote by $G_ \lambda , \ ( \lambda >0)$ the potential
(or resolvent)
operator of $X,$ which is given by $ G_ \lambda \psi (x) = E^x \int _0 ^ {+ \infty } e^{ - \lambda t } \psi (X(t)) dt,$
for any Borel function $\psi $ on $[0, + \infty ),$ and by $G^D_ \lambda $ the potential operator of  $\widetilde X$ on $D,$
and killed at $\tau_D (x),$ i.e. the hitting time of $\widetilde X$ to $0.$
We recall the following result from \cite{peng:2014}.

\begin{Theorem} \label{theoremzero}
Let $p^D$ denote the transition density of the diffusion $\widetilde X$ killed at zero; then, the Laplace transform of the FPT of the
DHJ  $X$ through $S$ is:
\begin{equation}
\widehat h _S ( \lambda , x ) = \frac {\widehat p(\lambda, x, S) } {\widehat p(\lambda, S, S) }  ,
\end{equation}
where $\widehat p$ is given by:
\begin{equation} \label{hatp}
\widehat p(\lambda, x, y)= \widehat p ^D (\lambda, x, y) + (1- \lambda G^D _ \lambda {\bf 1} (x) ) \cdot
\frac {\delta (y) + \beta \widehat p ^D (\lambda, \nu, y) } {\lambda (1+ \beta \langle \nu, G^D _ \lambda 1 \rangle ) }.
\end{equation}
\hfill $\Box$
\end{Theorem}
Let $\widetilde X$ be BM with drift $\mu  ,$ namely,
$\widetilde X(t) = \widetilde X_ \mu (t)  = x + \mu t + B_t ,$  then the transition density of
$\widetilde X$ killed at $0$ is (see \cite{cox:1965}, pg. 221):
\begin{equation}
p^D(t, x,y) = \frac 1 { \sqrt {2 \pi t }} \left [ e^{ - \left (y - \mu t -x \right )^2 /2t } - e^{-2 \mu x }
e^{ - \left (y-\mu t +x  \right ) ^2 /2t} \right ]  .
\end{equation}
By calculations analogous to those concerning the case $\mu =0$
(see e.g. \cite{revuz:1991}, pg. 82,  \cite{kar1:1975}, pg. 354,
\cite{kar2:1975}, pg. 288) we obtain
its Laplace transform:
\begin{equation} \label{hatpDWiener}
\widehat p ^D(\lambda, x, y)= \frac {e^{\mu (y-x)} } {\sqrt {2 \lambda + \mu ^2 } } \left [ e^{ - |y-x|  \sqrt {2 \lambda + \mu ^2 }
 } - e^{ -  |y+x| \sqrt {2 \lambda + \mu ^2 } } \right ] .
\end{equation}
Thus, if $\nu (dx) = g(x) dx,$
by straightforward calculations and using the expression of the Laplace transform
of the first hitting time of drifted BM to the barrier $S$ (see \cite{kar2:1975}), \eqref{hatp} becomes:
\begin{equation} \label{hatpWiener}
\widehat p(\lambda, x, y)= \widehat p ^D (\lambda, x, y) + e^{-x \left ( \sqrt {\mu ^2 + 2 \lambda} + \mu \right ) } \cdot
\frac {\delta (y) + \beta \int _ 0 ^\infty \widehat p ^D (\lambda, x, y) g(x) dx }
{\lambda +  \beta \left (1 - \int _0 ^\infty  e^{-x \left ( \sqrt {\mu ^2 + 2 \lambda} + \mu \right ) } g(x) dx \right )}
\end{equation}
where $\widehat p ^D$ is given by \eqref{hatpDWiener}. \par\noindent
Our main result is the following, where, for the sake of simplicity we limit ourselves to
consider the IFPT problem for the DHJ $X$ associated to BM with drift $\mu \le 0 ,$
in the special case when the starting point is $x=0.$

\begin{Theorem} \label{teoremauno}
Let  $\widetilde X(t)= \mu t + B_t , \ \mu \le 0 ,$ and let us consider the associated process $X$ with holding and jumps from $0,$
having holding rate $\beta >0$ and jump distribution $\nu;$
suppose that the FPT of $X$ over $S >0$
 has an assigned probability density  $f(t)= f_ \mu (t)$ and denote by
$ \widehat f _\mu(\lambda  ) = \int _0 ^ \infty e ^ {-\lambda t} f_\mu(t) dt , \ \lambda \ge 0 ,$
the Laplace transform of $f_\mu.$ Then, if
there exists a solution $(g, \beta)$ to the IFPT problem for $X,$  the Laplace transform
$ \widehat g (\lambda)$ of $g,$ and $\beta$  must satisfy the equation:

$$
\widehat f_\mu (\lambda   )
\Big [\left (1- e^{-2S \sqrt {2 \lambda + \mu ^2 }} \right ) \left ( \lambda
+ \beta \left (1-\widehat g \left ( \sqrt {2 \lambda + \mu ^2 } + \mu \right ) \right )
\right ) $$
$$+ \beta e^{ - S (2 \sqrt {\mu ^2 + 2 \lambda } \ + \mu ) }
\left ( \widehat g \left (- \sqrt {2 \lambda + \mu ^2 } \ \right )- \widehat g \left ( \sqrt {2 \lambda + \mu ^2 } \ \right ) \right ) \Big]$$
\begin{equation} \label{laplaceg}
=   \beta e^{-S (\sqrt {2 \lambda + \mu ^2 } \ - \mu ) }
\left ( \widehat g \left (- \sqrt {2 \lambda + \mu ^2 } \ \right )- \widehat g \left ( \sqrt {2 \lambda + \mu ^2 } \ \right ) \right )
\end{equation}
Moreover,  if  $\mu =0,$ and we require that the density $g$ is symmetric with respect to
$\frac S 2 ,$ then formula \eqref{laplaceg} can be explicited and we get:
\begin{equation} \label{laplaceesplicitag}
\widehat g (\lambda   ) = \frac {\widehat f _0  \left (\frac {\lambda ^2} 2 \right )
\left (\frac {\lambda ^2} 2  + \beta (1 + e^{-S \lambda }) \right ) } {\beta ( 1+  \widehat f _0  \left (\frac {\lambda ^2} 2 \right ) } \ .
\end{equation}
or
\begin{equation}
\widehat f _0 ( \lambda ) = \frac {\beta \widehat g ( \sqrt { 2 \lambda }) }
{\lambda + \beta (1+ e^{ -S \sqrt { 2 \lambda }} ) - \beta \widehat g ( \sqrt { 2 \lambda } ) }
\end{equation}
\end{Theorem}
{\it Proof.} \ Setting $\widehat f _ \mu (\lambda ) = \widehat h _S ( \lambda , 0 ),$ by Theorem \ref{theoremzero}, we have
\begin{equation}
\widehat f _ \mu (\lambda ) = \frac {\widehat p (\lambda, 0, S )} {\widehat p (\lambda, S, S ) }
\end{equation}
By using \eqref{hatpWiener} and calculating the various quantities, it is easy to see that:
$$ \widehat p (\lambda, 0, S )=  \frac {\beta } {\sqrt {2 \lambda + \mu ^2 } } \cdot
 \frac {e^{-S \left ( \sqrt {2 \lambda + \mu ^2 } \  - \mu \right ) }
\left ( \widehat g \left ( - \sqrt {2 \lambda + \mu ^2 } \ \right ) - \widehat g \left ( \sqrt {2 \lambda + \mu ^2 } \ \right ) \right )}
{\lambda + \beta \left (1- \widehat g \left ( \sqrt { 2 \lambda + \mu ^2 } + \mu \right ) \right ) } $$
and
$$ \widehat p (\lambda, S, S )= \frac {1-e^{-2S \sqrt {2 \lambda + \mu ^2 }} } {\sqrt {2 \lambda + \mu ^2 } } $$
$$+ \frac {\beta } {\sqrt {2 \lambda + \mu ^2 } } \cdot
 \frac {e^{-S \left ( 2 \sqrt {2 \lambda + \mu ^2 } \ + \mu \right )}
\left ( \widehat g \left ( - \sqrt {2 \lambda + \mu ^2 } \ \right ) - \widehat g \left ( \sqrt {2 \lambda + \mu ^2 } \ \right ) \right )}
{\lambda + \beta \left (1- \widehat g \left ( \sqrt { 2 \lambda + \mu ^2 } + \mu \right ) \right ) } $$
from which \eqref{laplaceg} follows. If $g$ is symmetric with respect to
$\frac S 2 ,$ we have:
\begin{equation} \label{hatgsimm}
\widehat g ( - \lambda )= e^{S \lambda } \widehat g (\lambda ).
\end{equation}
Then, equation \eqref{laplaceesplicitag}
easily follows, by taking $\mu =0$  in \eqref{laplaceg} and by using \eqref{hatgsimm}.
\hfill $\Box$
\bigskip
\begin{Remark}
{\rm If there is no holding at zero (i.e. $\beta = + \infty),$ by letting $\beta \rightarrow + \infty $ in
\eqref{laplaceesplicitag}, we obtain:
\begin{equation} \label{laplaceesplicitagnohold}
\widehat g (\lambda ) = \frac {\widehat f _0  \left (\frac {\lambda ^2} 2 \right ) (1 + e^{-S \lambda })  ) }
{1+  \widehat f _0  \left (\frac {\lambda ^2} 2 \right ) } .
\end{equation}
which coincides with equation (2.11) of \cite{abundo:saa13} . }
\end{Remark}
\begin{Remark}
{\rm By taking the first and second derivative with respect to $\lambda $ in \eqref{laplaceesplicitag} and calculating them
at $\lambda =0, $ we get
\begin{equation} \label{meantauS0}
E( \tau _S(0)) = \frac 1 \beta + 1 - 2 E( U^2)
\end{equation}
For $\beta \rightarrow + \infty, $ the above equation yields $E( \tau _S(0)) =  1 - 2 E( U^2)$
(cf. Remark 2.4 of \cite{abundo:saa13}).
}
\end{Remark}
\begin{Remark} {\rm A result analogous to that of Theorem \ref{teoremauno} can be also obtained for $x \neq 0;$ however, the
involved calculations are very heavy and cumbersome, so we have chosen not to  develop them. }
\end{Remark}

\begin{Remark} \label{remarkexistence}
{ \rm Once the pair $(\widehat g, \beta)$ has been found, such that it verifies \eqref{laplaceg}, or \eqref{laplaceesplicitag},
it may be that $\widehat g$ is not
the Laplace transform of
the density function of a random variable $U$ with support $(0,S).$ In this case, a solution to the IFPT problem does not exist.
This is the reason why Theorem \ref{teoremauno}
is formulated in a conditional form.
This kind of difficulty in showing the existence of a solution to an inverse FPT  problem is common to
other types of inverse problem (see e.g. \cite{abundo:saa14}, \cite{abundo:saa13}, \cite{abundo:stapro13},
\cite{abundo:stapro12}, \cite{abundo:saa06}, \cite{zuc:aap09});
as far as the present IFPT problem is concerned,
the difficulties are far stronger, because the relation between $\widehat g$ and $\widehat f$ is more complicated. For instance,
if $\mu =0$ and $f$ is the exponential density with parameter $1,$ then
the solution to the IFPT problem  with $S=1$ and $g$ symmetric in $(0,1)$
does not exist. In fact, suppose that the solution exists,
then by \eqref{laplaceesplicitag} it follows that
$$ \widehat g (\lambda ) = \frac {\lambda ^2 + 2 \beta  (1+ e^{- \lambda } ) } { \beta (4 + \lambda ^2 )}  ,$$
which is not the Laplace transform of a probability density in $(0,1),$ since  the third moment is negative. }
\end{Remark}
\par\noindent
Taking into account Remark \ref{remarkexistence}, we  will prove the existence of the density $g$ of $U$
for a class of FPT densities $f.$
For the sake of simplicity, we limit ourselves to the case when $\mu=0, \ x =0, \ S=1,$ and $g$ is required to be a function with support in $(0,1),$ which is symmetric with respect to the middle point $1/2 \ ;$
in fact, for $\mu \neq 0$ the calculations involved are far more complicated. \par\noindent
For an integer $k \ge 0,$ set $I_k(\lambda) = \int _{-1} ^1 e^ { - \lambda x } x^k dx;$ as easily seen,
$I _0 (\lambda )= 2 \sinh(\lambda ) / \lambda$ and the recursive relation
$ I _ {k} (\lambda ) = \frac { (-1)^k e^ \lambda - e ^ { - \lambda } } \lambda + \frac  k  \lambda I _ {k-1} (\lambda )$ allows to
calculate $I _k (\lambda),$ for every $k.$ \par
The following Proposition gives a sufficient condition, so that there exists the
solution $(g, \beta)$ to the IFPT problem for the process $X$ with holding and jumps from zero, associated to  $\widetilde X(t)=B_t,$ and the barrier $S=1.$

\begin{Proposition} \label{existenceproposition}
Let $\widetilde X(t)=B_t, $ and
suppose that the Laplace transform  of $f(t)$ has the form:
\begin{equation} \label{laplacedensityclass}
\widehat f ( \lambda )= \widehat f _ {2k}( \lambda ) =
\frac {(1 + \frac 1 {2k})e^{- \sqrt { \lambda /2 }} \left [ {\sqrt{2  / \lambda } } \sinh \left ( \sqrt { \lambda / 2} \right ) - I_{2k} \left (\sqrt { \lambda / 2 } \right ) \right ]}
{1+ e^{- \sqrt { 2 \lambda }} + \lambda / b - (1 + \frac 1 {2k} ) e^{- \sqrt { \lambda / 2 }} \left [ \sqrt{2  / \lambda }  \sinh \left (\sqrt { \lambda / 2} \right ) -
I_{2k} \left (\sqrt { \lambda / 2 } \right ) \right ] },
\end{equation}
for $b >0$ and some integer $k>0.$
Then, there exists the solution $(g, \beta)$ of the IFPT problem for $X,$  relative to the barrier $ S=1$ and
the FPT density $f,$ and
it results $\beta =b$ and:
\begin{equation}
 g(u)= g_{2k}(u)= \left (1 + \frac 1 { 2k} \right ) \left (1- (2u-1)^{2k} \right ), \ k \ge 0, \ u \in (0,1).
\end{equation}
\end{Proposition}
{\it Proof.}
A simple calculation shows that
$$\widehat g _ {2k} ( \lambda) = \left (1 + \frac 1 {2k} \right ) e^{ - \lambda /2} \left [ \frac 2 \lambda \sinh (\lambda /2) - I_{2k} ( \lambda /2) \right ] .$$
Since $g_{2k}$ is symmetric with respect to $S/2=1/2,$
the result follows by inserting $\widehat g _ {2k}$ into \eqref{laplaceesplicitag}.
\par \hfill  $\Box$ \par\noindent
Notice that, letting $b$ go to $+ \infty,$  \eqref{laplacedensityclass} becomes equation (2.17) of \cite{abundo:saa13}.
\bigskip

\begin{Remark} \label{remdopoesistenza}
A straightforward calculation shows that, if $U\in (0,1)$ has density $g_{2k},$ then $E(U^2)= \frac {4k+5 } { 6(2k+3)} .$
Then, by using \eqref{meantauS0},
 we obtain
that the FPT-distribution corresponding to $\widehat f _ {2k}$ has mean $E(\tau_1 (0))= \frac {2(k+2) } {3(2k+3) } + \frac 1 b \ .$
\end{Remark}

\section{Diffusions conjugated to Brownian motion}
In certain cases, a  one-dimensional diffusion $\widetilde X$ can be reduced to BM by a variable change; by using this approach,
we shall extend to a general one-dimensional DHJ $X,$  which is associated to $\widetilde X,$ the results obtained for Wiener process.
For $a >0,$ let  $J= [0, a]$ or $J= [0, a)$, with $a \le + \infty ,$ and  suppose that $\widetilde X$ is a time-homogeneous
diffusion in $J$
which is the solution of the stochastic differential
equation (SDE):

\begin{equation} \label{eqdiffu}
d \widetilde X(t) = \mu (\widetilde X(t)) dt + \sigma (\widetilde X(t)) d B_t \ , \  \widetilde X(0) = x \in J,
\end{equation}
where the coefficients $\mu (x)$ and $\sigma (x)$ are regular enough functions (see e.g. \cite{abundo:saa13}),
so that a unique strong solution exists.
We consider the following: \bigskip

\begin{Definition}
We say that $\widetilde X$
is conjugated to BM if there exists
an increasing  function $v :  J \longrightarrow \rm I\!R$ with $v(0)=0,$ such that: \par\noindent
(i) $v(x)$ is continuous for any $x \in  J $ and it is differentiable in the interior of $J;$ \par\noindent
(ii) $v^ {-1} (y)$ is differentiable in the interior of  $v(J)$ and it possesses the right derivative at \par\noindent
\ \ \ \ \ $y=0$ and the left derivative at $y=v(a),$ if $a < + \infty .$ \par\noindent
(iii)  $\widetilde X(t) = v^ {-1} (B_t + v(x)),$ for any $t \ge 0.$
\end{Definition}
\bigskip

\noindent Let us suppose that $\widetilde X$ is conjugated to BM, and consider the associated process $X$ with holding and
jumps from $0.$
Notice that holding and random reflection of $\widetilde X$ at zero corresponds to holding and random reflection of
$B_t + v(x)$ at zero;
moreover the first passage $\widetilde \tau _S (x)$ of $\widetilde X$ through the barrier $S,$ with $0<S<a, $ corresponds to the first passage of $B_t + v(x)$
through $v(S),$ and
$\widetilde \tau _S (x)= \widetilde \tau^B_{v(S)} (v(x)),$ where the superscript $B$ refers to BM. Furthermore, let $g(u)$ be
the density
of the  position $U\in (0,S)$ from
which, once $ X$ has hit $0,$ it starts afresh, after the holding time, and  let $q(y)$ be the corresponding density
of the position
$V= v (U)\in (0, v(S))$ from which, once $B_t + v(x)$ has hit $0,$ it starts afresh, after the holding time;
then, $q(y)=g(v^{-1} (y)) (v^{-1} )'(y),
\ y \in (0, v(S)).$  Thus,
if $\widetilde X$ is conjugated to BM via the function $v,$ then the solution $(g, \beta)$ to the IFPT problem
for the process
$X$ associated to  $\widetilde X,$ relative to the FPT density $f$
and the barrier $S,$ can be written in terms of the solution $(q, \beta)$ to the IFPT problem for the process associated to
$v(x)+ B_t,$
relative to the FPT density $f$
and the
barrier $v (S),$  by using that $g(x)= q(v(x)) v'(x).$  As easily seen, if $x=0, \ \widehat q$ is obtained by \eqref{laplaceg} with $\mu =0,$
with $\widehat q$ in place of $\widehat g$ and $v(S)$ in place of $S.$ \par\noindent

\section{A few examples}
\noindent {\bf Example 1} \ (when  $g$ is the uniform density in $(0,S)$)\par\noindent
Let $\widetilde X(t)= B_t$ and $S >0, \ b >0$ and let
$$\widehat f(\lambda )= \frac {b (1-e^{-S \sqrt {2 \lambda }} )}
{ S \sqrt { 2 \lambda } [ \lambda +   b (1 + e^{-S \sqrt { 2 \lambda } } )] - b (1 - e^{- S \sqrt {2 \lambda }} )  }  \ .$$
Then, the solution to the IFPT problem for $X$  relative to $S$ is the pair $(g, \beta )$ with $\beta =b$ and $g(u)= \frac 1 S {\bf 1} _ {(0,S)} (u),$
i.e. the uniform density in $(0,S);$
this is easily
obtained by searching for a solution which is symmetric with respect to $S/2;$ by using \eqref{laplaceesplicitag} with $\beta = b,$  we get
$\widehat g(\lambda )= \frac {1-e^{-S \lambda } } {S \lambda  } $ which is the
Laplace transform of $g(u) = \frac 1 S {\bf 1} _ {(0,S)}(u) .$
In the case $S=1,$ we can obtain the same result by letting $k$ go to infinity in $\widehat f_ {2k}$ and $g_{2k}$ of
Proposition \ref{existenceproposition};
moreover, the mean of the FPT-distribution corresponding to $\widehat f$  is $1/3 + 1/ b,$
as it also follows by calculating $\lim _ {k \rightarrow \infty}  \left (  \frac {2(k+2) } {3(2k+3) } + \frac 1 b \right )$
(see Remark \ref{remdopoesistenza} ). For $b \rightarrow + \infty$ (no holding at $0),$ one obtains the function
$ \widehat f( \lambda )$ of Example 3 of \cite{abundo:saa13}.
\bigskip

\noindent{\bf Example 2} \
Let $\widetilde X(t)= B_t$ and $S >0, \ b >0$ and let
$$ \widehat f( \lambda ) = \frac {b \pi ^2 (1+ e^{-S \sqrt { 2 \lambda }} ) }
{(4 \lambda S^2 + 2 \pi ^2 )( \lambda + b (1+ e^{-S \sqrt { 2 \lambda }} )) - b \pi ^2 (1+ e^{-S \sqrt { 2 \lambda }}  )};
$$
then the solution to the IFPT problem for $X$  relative to
$S$ is the pair $(g, \beta)$ with $\beta =b$ and $g(u)= \frac { \pi } {2S} \sin \left ( \frac \pi S  u \right ), \ u \in (0,S) .$ \par\noindent
In fact, we search for
a solution which is symmetric with respect to $S/2;$ by using \eqref{laplaceesplicitag} with $\beta =b,$ it follows that
$\widehat g (\lambda ) = \frac {\pi ^2} 2 \ \frac {\left (1+e^{ - \lambda S}\right ) } {\lambda ^2 S^2 + \pi ^2 },$
which is indeed the Laplace transform of the function $g(u)$ above. If there is no holding at $0$ (i.e. $b = + \infty ),$
$ \widehat f( \lambda )$ becomes that of Example 1 of \cite{abundo:saa13}.
\bigskip

\noindent {\bf Example 3} (when  $g$ is a Beta density) \par\noindent
Let $\widetilde X(t)= B_t$ and  $S>0, \ b >0$ and let
$$\widehat f(\lambda )= \frac {6(e^{-S\sqrt {2 \lambda}} (S \sqrt {2 \lambda } +2 ) + S \sqrt { 2 \lambda } -2 ) }
{S^3 \lambda ^3 (1+e^{-S \sqrt { 2 \lambda }} ) - 6(e^{-S\sqrt {2 \lambda}} (S \sqrt {2 \lambda } +2 ) + S \sqrt { 2 \lambda } -2 ) + \lambda / b } \ .$$
Then, the solution to the IFPT problem for $X$  relative to $S$ is the pair $(g, \beta)$ with $\beta = b$ and $g(u)=\frac {6} { S^3} \ u(S-u), \ u \in (0,S) .$
In fact, by using \eqref{laplaceesplicitag} with $\beta =b,$  it follows  that
$\widehat g (\lambda ) = \frac 6 {S^3 \lambda ^3 }  [ e^{ - S \lambda } (S \lambda +2 ) + S \lambda -2],$
which is indeed the Laplace transform of $g(u)=\frac {6} { S^3} \ u(S-u) .$
Notice that, for $S=1, \ g$ is the density $g_{2k}$ of Proposition \ref{existenceproposition}, for $k=1.$
For $b \rightarrow + \infty$ (no holding at $0),$ one obtains the function
$ \widehat f( \lambda )$ of Example 4 of \cite{abundo:saa13}.

\bigskip

\noindent {\bf Example 4} (when $g$ is the triangular density in $(0,1)$) \par\noindent
Let $\widetilde X(t)= B_t$ and $S=1, \ b >0$ and let
$$\widehat f(\lambda )= \frac {2 b (1- e^{ \sqrt { \lambda /2}}  )^2 }
{ \lambda [ \lambda  + b (1 + e^{- \sqrt { 2 \lambda }} )] - 2 b (1 - e^{- \sqrt { \lambda / 2 } } )^2
} $$
Then, the solution to the IFPT problem for $X$  relative to $S$ is the pair $(g, \beta )$ where $\beta =b$ and $g$ is the triangular density in $(0,1):$
$$ g(u) =
\begin{cases}
4u, \ u \in (0, \frac 1 2 ] \\ 4(1-u) , \ u \in ( \frac 1 2 , 1)
\end{cases} \ .
$$
In fact, by using \eqref{laplaceesplicitag} with $\beta =b,$ it follows  that $\widehat g (\lambda ) =
\frac 4 {\lambda ^2 }  (1 - e ^{ - \lambda /2 }) ^2,$  which is indeed the Laplace transform of the function $g(u)$ above.
For $b \rightarrow + \infty$ (no holding at $0),$ one obtains the function
$ \widehat f( \lambda )$ of Example 5 of \cite{abundo:saa13}.

\bigskip

\noindent {\bf Example 5.} \par\noindent
A class of diffusions  conjugated to BM is given by processes $\widetilde X(t)$ which are solutions of SDEs such as
\begin{equation} \label{conjdiffu}
 d\widetilde X(t) = \frac 1 2 \sigma (\widetilde X(t)) \sigma ' (\widetilde X(t)) dt + \sigma (\widetilde X(t)) dB_t, \  \widetilde X(0)=x
 \end{equation}
with $\sigma (\cdot) \ge 0.$
Indeed, if the integral
$v(z) \doteq  \int _{x} ^ z \frac {1}  {\sigma (r) }dr $
is convergent for every $z ,$
by It${\rm \hat o}$'s  formula, we obtain that $\widetilde X(t)= v^ {-1} (B_t +v(x)).$
Then, by using the results of Section 2, examples of solutions to IFPT problems with $x=0$ can be easily derived
from Examples 1 to 4, with regard to the process $X(t)$
with holding and jumps from $0,$ associated to $\widetilde X$ which is driven by the SDE \eqref{conjdiffu},
for some choice of $\sigma (\cdot );$ in fact, it suffices to replace
$S$ with $v(S).$ \par\noindent
For instance, diffusions $\widetilde X$ of this kind are the well-known
Feller process (also known as the Cox-Ingersoll-Ross (CIR) model), and the
Wright \& Fisher-like process (see Examples (i) and (ii) of \cite{abundo:saa13}).
\bigskip

\noindent{\bf Example 6} (Ornstein-Uhlenbeck process)  \par\noindent
Let  $\widetilde X(t)$ be the solution of the SDE:
$$ d\widetilde X(t) = - \mu \widetilde X (t) dt + \sigma dB_t, \ \widetilde X(0)= x ,$$
where $\mu, \sigma$ are positive constants. By using a time--change, the explicit solution assumes the form
$ \widetilde X(t)= e^{- \mu t } \left ( x + B( \rho(t)) \right ),$  where
$ \rho (t) =  \frac { \sigma ^2 } { 2 \mu } \left ( e ^ {2 \mu t} -1 \right).$
If $S(t)$ is a moving barrier, the FPT of $\widetilde X (t)$ over $S(t)$ is
$ \widetilde \tau _{S(t)}   =  \inf \{ t>0: x+ B (\rho (t)) \ge e^ {\mu t } S(t) \} $
and so
$ \rho ( \widetilde \tau_{S(t)}) =  \inf \{ u >0: x + B_u \ge  \widetilde S (u) \} ,$
where $\widetilde S (u)  = e^ {\mu \rho ^{-1} (u) } S(\rho ^{-1} (u)) .$
Therefore, if  $S(t)= S_0e^{ - \mu t},$
the IFPT problem for the associated DHJ $X(t),$  and
relative to the moving barrier $S(t)$ and the FPT distribution $F,$
is reduced to the IFPT problem for
DHJ associated to BM, starting from $x$ and relative to the  constant barrier $S_0$ and the FPT distribution $\widetilde F = F \circ \rho ^ {-1} .$
Thus, if $x=0,$ explicit examples for the Ornstein-Uhlenbeck process and the exponential barrier $S(t)= S_0e^{ - \mu t},$ can be
easily derived from Examples 1 to 4.
\bigskip

\noindent {\bf Example 7} (Geometric Brownian motion)  \par\noindent
Let $\widetilde X(t)$ be the solution of the SDE:
$$d\widetilde X(t) = r \widetilde X(t) dt + \sigma \widetilde X(t) d B_t , \ \widetilde X(0)= x >0  , $$
where $ r$ and $\sigma$ are positive constant. This is a well-known equation in the framework of Mathematical Finance, since
it describes the time evolution of a stock price $\widetilde X.$ The explicit solution is
$\widetilde X(t)= x e ^ { \mu t } e ^ { \sigma B_t} ,$
where $\mu = r - \sigma ^2 /2 .$
Let us consider the moving barrier $S(t)= e ^ {\sigma S + \mu ' t };$
then, the IFPT problem for the associated DHJ $X(t),$ relative to $S(t)$ and  the FPT distribution $F,$ is reduced to
the IFPT problem for DHJ associated to BM with drift $( \mu - \mu ')/ \sigma ,$ starting from $\frac {\ln x } \sigma $ and relative to the constant  boundary $S ,$
and the same FPT distribution $F.$


\begin{thebibliography}{spc}

\bibitem [1] {abundo:saa14}
Abundo, M. 2014.\newblock
One-dimensional reflected diffusions with two boundaries and an inverse first-hitting problem.
\newblock{\it Stochastic Anal. Appl.} 32: 4, 975–991, DOI: 10.1080/07362994.2014.959595

\bibitem [2] {abundo:saa13}
Abundo, M. 2013.\newblock
Solving an inverse first-passage-time problem for Wiener process subject to
      random jumps from a boundary.
\newblock{\it Stochastic Anal. Appl.} 31: 4, 695-707.


\bibitem [3] {abundo:stapro13}
Abundo, M. 2013.\newblock
The double-barrier inverse first-passage problem for Wiener process
with random starting point.
\newblock{\it Statist. Probab. Lett.} 83: 168--176.


\bibitem [4] {abundo:stapro12}
Abundo, M. 2012.\newblock
An inverse first-passage problem for one-dimensional diffusions with random
starting point.
\newblock{\it Statist. Probab. Lett. } 82(1): 7--14.

\bibitem [5] {abundo:saa11}
Abundo, M. 2011.\newblock
First passage problems for one-dimensional diffusions with random jumps  from a boundary.
\newblock{\it Stochastic Anal. Appl.} 29(1): 121--145.

\bibitem [6] {abundo:saa06}
Abundo, M. 2006. \newblock
Limit at zero of the first-passage time density and the inverse problem for
one-dimensional diffusions.
\newblock{\it Stochastic Anal. Appl.} 24: 1119--1145.





\bibitem[7]{abundo:pms97}
Abundo, M. 1997. \newblock
On some properties of one-dimensional diffusion processes on an interval.
\newblock   {\it Prob.  Math.  Statis.} 17(2): 235--268.

\bibitem[8]{alb:2007}
Albano, G., Giorno, V., Nobile, A.G., and Ricciardi, L.M. 2007. \newblock
A Wiener-type neuronal model in the presence of exponential refractoriness.
\newblock   {\it BioSystem} 88: 202--215.

\bibitem[9]{benari:spa08}
Ben-Ari, I., Pinsky, R.G. 2009. \newblock
Ergodic behavior of diffusions with random jumps from the boundary. \newblock  {\it Stoch. Processes  Appl.} 119(3): 864--881.

\bibitem[10]
{cox:1965}
Cox, D.R. and Miller, H.D., 1965.\newblock
{\it The Theory of Stochastic Processes.} \newblock
Methuen $\&$ Co LTD, London.



\bibitem[11]
{dicresc:que03}
Di Crescenzo, A., Giorno, V., Nobile, A.G. and Ricciardi, L.M. 2003. \newblock
On the M/M/1 queue with catastrophes and its continuous approximation. \newblock {\it Queueing System} 43: 329--347.


\bibitem [12] {feller:1954}
Feller, W. 1954.\newblock
Diffusion processes in one dimension.
\newblock{\it Trans. Amer. Math. Soc.} 17: 1--31.




\bibitem[13]
{jackson:stapro09}
Jackson, K, Kreinin, A, and Zhang, W. 2009.\newblock
Randomization in the first hitting problem. \newblock {\it Statist. Probab. Lett.} 79: 2422--2428.



\bibitem[14]
{kar1:1975}
Karlin, S. and Taylor, H.M. 1975.\newblock
{\it A first course in stochastic processes.} \newblock
Academic Press, New York.

\bibitem[15]
{kar2:1975}
Karlin, S. and Taylor, H.M. 1981.\newblock
{\it A second course in stochastic processes.} \newblock
Academic Press, New York.

\bibitem[16]{lansky:jtb94}
Lanska, V., Lansky, P. and Smiths, C.E. 1994. \newblock
Synaptic transmission in a diffusion model for neural activity.
\newblock   {\it J. Theor. Biol.} 166: 393--406.

\bibitem[17]{lansky:bio89}
Lansky, P. and Smiths, C.E. 1989. \newblock The effect of a random
initial value in neural first- passage-time models.
\newblock   {\it Math. Biosci.} 93(2): 191–-215.


\bibitem[18]{peng:2014}
Peng, J. 2014. \newblock
A note on the first passage time of diffusions with holding and jumping boujndary.
\newblock   {\it Statist. Probab. Lett.} 93: 58--64.


\bibitem[19]{peng:2013}
Peng, J. and Li, W. 2013. \newblock
Diffusions with holding and jumping boundary.
\newblock   {\it Sci.China Math.} 1: 161--176.

\bibitem[20]{peng:sta09}
Peng, J. and Zaiming Liu. 2009. \newblock
On a class of mathematical ecosystems with random jumps.
\newblock   {\it Statist. Probab. Lett.} 79(5): 630--636.

\bibitem[21]
{revuz:1991}
Revuz, D. and Yor, M., 1991.  \newblock
{\it Continous martingales and Brownian motion.} \newblock
Springer-Verlag, Berlin
Heidelberg.







\bibitem [22] {zuc:aap09}
Zucca, C. and Sacerdote, L. 2009. \newblock
On the inverse first-passage-time problem for a Wiener process. \newblock
{\it Ann. Appl. Probab.} 19 (4): 1319–-1346.


\end{thebibliography}
\end{document}